\magnification 1200
\def\R{{\rm I\kern-0.2em R\kern0.2em \kern-0.2em}}
\def\N{{\rm I\kern-0.2em N\kern0.2em \kern-0.2em}}
\def\P{{\rm I\kern-0.2em P\kern0.2em \kern-0.2em}}
\def\B{{\rm I\kern-0.2em B\kern0.2em \kern-0.2em}}
\def\Z{{\rm I\kern-0.2em Z\kern0.2em \kern-0.2em}}
\def\C{{\bf \rm C}\kern-.4em {\vrule height1.4ex width.08em depth-.04ex}\;}

\def\D{{\Delta}}

\def\z{{\zeta}}

\font\ninerm=cmr8
\  
\vskip 20mm
\centerline {\bf ANALYTICITY ON TRANSLATES OF A JORDAN CURVE}
\vskip 4mm
\centerline{Josip Globevnik}
\vskip 4mm
{\noindent \ninerm ABSTRACT\ \   Let $\Omega $ be a domain in $\C $ which
is symmetric with respect to 
the real axis and whose boundary is a real analytic simple closed curve.
Translate $\overline\Omega $ vertically to get 
$K = \cup\{ \overline\Omega +it,\ -r\leq t\leq r\}$ where $r>0$ is such that 
$(\overline\Omega -ir)\cap(\overline\Omega +ir)=
\emptyset $. We prove that if $f$ is a continuous function on $K$ 
such that 
for each $t,\ -r\leq t\leq r$, the function $f|(b\Omega+it)$ has a 
continuous extension to 
$\overline\Omega +it$ which is holomorphic on $\Omega +it$ 
then $f$ is holomorphic on $\hbox{\rm Int}K$.} 
\vskip 4mm

\bf 1.\ Introduction 
\vskip 2mm
\rm
Write $\D (a,r)=\{\z\in\C\colon\ |\z -a|<r\}$ and let
$\D =\D (0,1)$. Translate $b\D $ vertically 
to get the strip 
$S=b\D +i\R = \{ \z\in \C :\ -1\leq \Re \z \leq 1\}$. Let $f$ be a continuous 
function on $S$ such that for each $t\in\R $ the function $f|(b\D +it)$ has a 
continuous extension to 
$\overline \D +it$ which is holomorphic on $\D +it$. Must $f$ be holomorphic on  
$\hbox{Int}S$ [G2]? 

A positive answer was obtained for real analytic functions by M.\ Agranovsky and 
the author [AG] and independently by L.\ Ehrenpreis [E]  
and for  
continuous functions by A.\ Tumanov [T1].

To answer the question above one passes in both [AG] and [T1] to 
an associated problem in $\C ^2$. 
In [AG] the authors use semi-quadrics 
$$
\Lambda _{a,\rho} = \{ (z,w)\colon\ (z-a)(w-\overline a)=\rho^2,\ 
0<|z-a|<\rho\}
$$
which are attached to $\Sigma = \{(z ,\overline z )\colon\ z\in\C\}$ along the circles
$\{(z,\overline z)\colon\ z\in b\D (a,\rho )\}$. They use the property of $\Lambda _{a,\rho}$ 
that a continuous function on $b\D (a,\rho )$ extends holomorphically through 
if and only if the function $F$ defined on $b\Lambda _{a,\rho} = \{ (z ,\overline z )\colon\ 
z \in b\D (a,\rho)\}$ by $F(z,\overline z)= f(z)\ (z\in b\D (a,\rho ))$ has a bounded 
continuous extension 
to $\Lambda _{a,\rho}\cup b\Lambda _{a,\rho}$ which is holomorphic on $\Lambda _{a,\rho}$.
Thus, 
when studying holomorphic extensions of a function $f$ from circles in the plane one 
defines 
$F(z,\overline z)= f(z)$ in a region in $\Sigma$ and then studies bounded holomorphic 
extensions of $F$ from $b\Lambda _{a,\rho}$ through $\Lambda _{a,\rho}$ [AG, G3].

Tumanov [T1] passes to a problem in $\C^2$ by adding  
an extra variable to make the translates of the disc pairwise disjoint
and then, on the union of these discs, the (smooth) manifold
$$
N_1 = \{(\z+it,\z )\colon\ \z\in\overline\D,\ t\in\R \},
$$
he defines a continuous function $F$ by letting, for each $t\in\R$, the function  
$\z\mapsto F(\z +it,\z)\ (\z\in\overline\D)$ be 
the holomorphic extension of $\z\mapsto f(\z+it)\ (\z\in b\D)$ 
through $\D$. 
He then observes that the symmetry
$$
F(\z +it,\z)= F(\z+it, 1/\z ) \ \ (\z\in b\D , t\in \R)
\eqno (1.1)
$$
which follows from the fact that $F(\z +it, \z))=f(\z+it) \ (\z\in b\D,\ t\in\R)$, 
makes possible to extend $F$ continuously to a new geometric object, the (smooth) manifold
$$
N_2 = \{ (\z +it,1/\z )\colon\ \z\in\overline\D\setminus \{ 0\},\ t \in \R),
$$
by using the equality (1.1) for $\z\in\overline\D\setminus \{ 0\},\ t\in\R $ as
a definition. 
Thus one gets a continuous CR function $F$ on $N_1\cup N_2$, the union of two
manifolds with common boundary
$$
N_1\cap N_2 = bN_1 =bN_2 = \{ (\z+it, \z)\colon \z\in b\D ,\  t\in \R\}.
$$
Tumanov then uses methods of CR theory to show that $F$ does not depend on the second 
variable which means 
that $f$ is holomorphic on $\hbox{Int}S$. He also discovers that this is 
actually a finite strip problem. 

Very recently Tumanov [T2] studied a similar problem for a family of circles
with centers sliding along a smooth curve and 
with smoothly changing radii. He used semi-quadrics. He obtained the result 
by using a classical argument of H.\ Lewy about 
holomorphic extensions of CR functions. In particular, he found a very simple
proof of the theorem on the strip. 

In the present paper we generalize the
result of Tumanov [T1] from vertical translates of circles  
to vertical translates of real-analytic simple closed curves which 
are symmetric with respect
to the real axis. 

\vskip 4mm
\bf 2.\ The main results \rm
\vskip 2mm
Our first result is about analyticity on 
vertical translates of curves which are symmetric with respect to the real axis.
\vskip 2mm
\noindent\bf Theorem 2.1\ \it Let $\Omega $ be a domain in $\C $ which
is symmetric with respect to 
the real axis and whose boundary is a real analytic simple closed curve.
Let  $r>0$ and let 
$f$ be a continuous function on $K = \cup\{ b\Omega +it,\ -r\leq t\leq r\}$
such that 
for each $t,\ -r\leq t\leq r$, the function $f|(b\Omega+it)$ has a 
continuous extension to 
$\overline\Omega +it$ which is holomorphic on $\Omega +it$. Suppose 
that $(\overline\Omega -ir)\cap(\overline\Omega +ir)=
\emptyset $. Then $f$ is holomorphic on $\hbox{\rm Int}K$. \rm
\vskip 2mm
\noindent Note that our assumptions imply that $K = \cup\{ \overline\Omega
+it,\ -r\leq t\leq r\}$.

We will deduce Theorem 2.1 from a more general result below which involves vertical 
translates of general domains 
 and 
their images under conjugation.

Let $D$ be a domain in $\C$ bounded by a real-analytic simple closed curve. 
Let $S$ be the vertical strip defined by 
  $S=\cup\{bD+it,\ t\in\R\} = \{ \z\in\C\colon\ \alpha\leq\Re \z\leq \beta\}$. 
 Write
$D^\ast = \{\overline\z\colon \ \z\in D\}$. Obviously,  
$\cup\{b{D^\ast}+is,\ s\in\R\}= S$.
\vskip 2mm
\noindent\bf
Theorem 2.2\ \it Let $\lambda\colon [\alpha,\beta]\rightarrow \R$ be a
continuous function and let 
$a, b, c, d$ be real numbers such that $\overline D+ia,\ \overline {D^\ast}+ic$ are 
both contained in 
$\{ t+is\colon\ s<\lambda (t),\ \alpha\leq t\leq\beta \}$ and such that $\overline D+ib,\ 
\overline {D^\ast}+id$ are both contained in 
$\{ t+is\colon\ s>\lambda (t),\ \alpha\leq t\leq\beta \}$. Let 
$$
Q_1 =\cup \{ bD+it\colon\ a\leq t\leq b\},
\ Q_2 =\cup \{ b{D^\ast}+is\colon\ c\leq s\leq d\}
$$
and let $f$ be a continuous function on $Q_1\cup Q_2$ such that
$$
\left.\eqalign {&\hbox{for each\ } t, \ \ a\leq t\leq b,\ \hbox{the function\ }
f|(bD+it)\hbox{\ has a conti-}\cr &\hbox{nuous extension to\ } 
\overline D+it \hbox{\ which is holomorphic on\ }D+it,\cr}
\right\} \eqno (2.1)
$$
$$
\left.\eqalign {\hbox{for each\ } s, \ \ c\leq s\leq d,\ \hbox{the function\ }
f|(bD^\ast+is)\hbox{\ has a conti-}\cr \hbox{ nuous extension to\ } 
\overline {D^\ast}+is \hbox{\ which is holomorphic on\ }D^\ast +is.\cr}
\right\} \eqno (2.2)
$$
Then the function $f$ is holomorphic on $\hbox{\rm Int}Q_1\cup \hbox{\rm Int}Q_2$.\rm 
\vskip 2mm
\noindent Our assumptions about $a, b, c, d$ mean that 
$\overline D+ia,\ \overline{D^\ast} 
+ ic$ both lie below the curve $\ell = 
\{ t+i\lambda (t)\colon\ \alpha\leq t\leq\beta\}$ and
$\overline D+ib,\ \overline{D^\ast} 
+ id$ both lie above the curve $\ell$. Note that this implies that 
$Q_1 =\cup \{ \overline D+it\colon\ a\leq t\leq b\}$ and 
\ $Q_2 =\cup \{ \overline {D^\ast}+is\colon\ c\leq s\leq d\}$.
Theorem 2.1 follows from
Theorem 2.2 by putting 
$\Omega = D = D^\ast$, $a=c=-r, \ b=d=r$ and $\lambda\equiv 0$, 
that is, $\ell =[\alpha ,\beta]$. 
\vskip 4mm
\bf 3.\ From circles to general curves

\rm\vskip 2mm In this section we describe the idea how to pass from circles to general curves. 
Let $\Omega $ be a domain bounded by a real-analytic simple closed curve 
which is symmetric with respect to the real axis. With no loss of generality assume 
that $\Omega$
contains the origin. Let $f$ be a continuous function on $\cup\{b\Omega+it\colon\ t\in\R\}$ 
such that for each $t\in\R$, the function $\z\mapsto f(\z +it)\ (\z\in b\Omega)$ 
has a continuous 
extension to $\overline\Omega $ which is holomorphic on $\Omega $. 

Semi-quadrics are related to circles so they cannot be used 
to study the analyticity of functions 
on a family of translates
of a given curve that is not a circle. We look again at the way how Tumanov [T1] adds the 
extra variable in the case of the circles. An important point in his setting  
is that on $b\D$ the conjugation $z\mapsto\overline z$ extends to the map $z\mapsto 1/z$ 
which carries $\Delta
\setminus \{ 0\}$ biholomorphically onto $\C\setminus\overline\D$. This is not the
case for general curves 
so for domains $\Omega $ more general than a disc it seems difficult to work with the manifold
$\{ (\z+it, \z)\colon\ \z\in
\overline\Omega, t\in\R\}$ used in [T1] when $\Omega $ is a disc. However, since 
the reflection (1.1) takes place
only in 
the second variable the  
idea is to replace the manifold  $\{ (\z+it, \z)\colon\ \z\in
\overline\Omega, t\in\R\}$   with a manifold that is attached to the cylinder $
\{(z,w)\colon \ |w|=1\}$. To do this we take the conformal 
map $\Phi \colon\Omega\mapsto \D$ 
that satisfies  
$\Phi (\overline\z ) =\overline{\Phi (\z )} \ (\z\in\Omega ),\ \Phi (0)=0$, 
notice that $\Phi $ 
extends to a diffeomorphism $\Phi $ from $\overline\Omega $ to $\overline \D$ 
and define  the smooth manifold $N_1$ by 
$$
N_1=\{ (\z+it,\Phi (\z ))\colon \ \z\in\overline\Omega,\ t\in\R\}.
$$
We define a continuous function $F$ on $N_1$ by letting, for each 
$t\in\R$, the function $\z\mapsto F(\z+it,\Phi (\z ))\ (\z\in\overline\Omega)$ 
be  the 
holomorphic extension of $\z\mapsto f(\z+it)\ (\z\in b\Omega)$ through $\Omega $. 

If $\z \in b\Omega $ then $\overline\z\in b\Omega $ and if $t\in\R$ then 
$\z+it = \overline\z+is$ where $s=t+(\z-\overline\z )/i \in\R$ so $F(\z+it,\Phi (\z))=f(
\z +it)=f(\overline\z +is)=F(\overline\z+is,\Phi (\overline\z ))
= F(\overline \z +is,\overline{\Phi (\z )})
= F(\overline\z +is,1/\Phi (\z )) = F(\z+it,1/\Phi (\z ))$ so
$$
F(\z+it,\Phi (\z ))= F(\z+it,1/\Phi (\z ))\ \ \ (\z\in b\Omega,\ t\in\R)
\eqno (3.1)
$$
which makes possible to extend $F$ continuously to a new geometric object, the smooth manifold
$$
N_2=\{(\z+it,1/\Phi (\z ))\colon\ \z\in\overline\Omega \setminus\{ 0\},\ t\in\R\},
$$
by using the equality (3.1) for $\z\in\overline\Omega \setminus\{ 0\}$ as a
definition. Thus we 
get a continuous CR function $F$ on $N_1\cup N_2$, the union of two 
manifolds with the common boundary 
$$
N_1\cap N_2=bN_1=bN_2=\{ (\z+it,\Phi (\z ))\colon\ \z\in b\Omega ,\ t\in\R\}.
$$
We then show that the classical argument of H.\ Lewy which Tumanov 
used with semiquadrics works also in the present situation. This helps us to prove that $F$ 
depends only on the first variable which implies that $f$ is holomorphic. 

In fact, our main result, Theorem 2.2, is somewhat more general than the one just 
described. Its proof, although technically a bit complicated, uses essentially the idea above. 
\vskip 4mm
\noindent \bf 4.\ The manifold $N$\rm
\vskip 2mm
We now begin with the proof of Theorem 2.2. With no loss of generality we assume
that $0\in D$ and 
that the imaginary axis intersects $bD$ transversely. 

Let $\Phi\colon\ D\rightarrow\D $ be a conformal map such that $\Phi (0)=0$.
Since $bD$ is 
real-analytic the map $\Phi $ extends to a biholomorphic map from a neighbourhood 
of $\overline D$ to a neighbourhood 
of $\overline\D$. Define $\Psi\colon\ D^\ast \rightarrow \D $ by
$$
\Psi (\z ) =\overline{\Phi (\overline\z )}\ \ (\z\in D^\ast ).
$$
The map $\Psi $ maps $D^\ast $ conformally onto $\D $ and extends
to a biholomorphic map from a neighbourhood of $\overline{ D^\ast }$ to a neighbourhood 
of $\overline\D$.

Define
$$\eqalign{
&N_1=\{(\xi +it, \Phi (\xi))\colon\ \xi\in\overline D,\ t\in\R \},\cr 
&N_2=\{(\z +is, 1/\Psi (\z))\colon\ \z\in\overline {D^\ast}\setminus\{  0\},\ s\in\R \}\cr}
$$
and set $N=N_1\cup N_2$. Write
$$
\Phi ^{-1}(w)=p(w)+iq(w)\ \ (w\in\overline\D )
$$
where $p$ and $q$ are real functions. Then $N_1 = \{ (p(w)+iq(w)+it,w)\colon\ w\in
\overline\D,\ t\in\R\}$ $= \{ (p(w)+it,w)\colon\ w\in\overline\D,\ t\in\R\} = 
\{ (p(w), w)\colon\ w\in\overline\D\} + \R (i,0)$. If $w=1/\Psi (\z )$ then 
$\overline{\Phi (\overline\z)}= 
\Psi(\z )=1/w$ so $\z=\overline{\Phi^{-1}(1/\overline w)}=
p(1/\overline w)-iq(1/\overline w)$ which implies that
$N_2=\{(p(1/\overline w)-iq(1/\overline w)+is, w)\colon \ w\in \C\setminus\Delta ,\ s\in\R\} =
\{ (p(1/\overline w)+is,w)\colon\ w\in\C\setminus \D,\ s\in\R\} =
\{ p(1/\overline w),w)\colon\ w\in\C\setminus \D\}+\R(i,0)$.
Define
$$
\theta (w)=\left\{\eqalign{&p(w)\ \ \ \ \ \ \ (w\in\overline\D )\cr
&p(1/\overline w )\ \ \ \ (w\in \C\setminus \D ).\cr}\right. 
$$
The function $\theta $ is well defined and continuous on
$\C$ since \ $w=1/\overline w\ \ (w\in b\D)$. 
The function $\theta $ is invariant with respect to $w\rightarrow 1/\overline w$, 
the reflection 
across $b\D $. Note that $\theta|\overline \D$ is smooth as it extends to a 
harmonic function in a neighbourhood of $\overline\D$. Similarly, $\theta|(\C\setminus\D )$ 
is smooth. We have
$$
N= \{ (\theta (w)+it,w)\colon\ w\in\C,  t\in\R\} 
\eqno (4.1)
$$
which shows that we obtain $N$ by taking the graph $\{(\theta (w),w)\colon\ w\in\C\}$ 
of $\theta $ in $\R\times\C =
\R\times\{ 0\}\times\C$ and then making the union of all translates of this graph in the
extra perpendicular direction
$(i,0)$, that is, 
$$
N=\{ (\theta (w),w)\colon\ w\in\C\} + \R (i,0).
$$
Since 
$$
\eqalign{
&N_1= \{ (\theta (w),w)\colon\ w\in\overline\D\} +\R (i,0)\cr
&N_2= \{ (\theta (w),w)\colon\ w\in\C\setminus\D\} +\R (i,0)\cr}
\eqno (4.2)
$$
we see that $N$ is the union of manifolds $N_1 $ and $N_2$ with 
boundary which meet along the common boundary
$$
N_1\cap N_2 = bN_1=bN_2=\{(\theta (w),w)\colon\ w\in b\D \} + \R (i,0).
$$

The complement of the graph of $\theta $ in $\R\times\C$ has two components:
$\{ (t,w)\colon\ t>\theta (w),\ w\in\C\} $ and 
$\{ (t,w)\colon\ t<\theta (w),\ w\in\C\} $, which, by (4.1) implies that $\C^2\setminus N$ 
has two components
$$
\eqalign{
&P_1=\{(t+is,w)\colon \ t>\theta (w),\ s\in \R, w\in\C\} \cr
&P_2=\{(t+is,w)\colon \ t<\theta (w),\ s\in \R, w\in\C\}. \cr}
$$
\vskip 4mm
\bf 5.\  Intersecting $N$ with complex lines \rm
\vskip 2mm
We will apply the reasoning of H.\ Lewy about holomorphic extensions of CR 
functions. To this end, we 
look first at the intersections of $N$ with complex lines $L(z)=\{ (z,w)\colon\ w\in\C \}$.
We 
shall use the fact that since $bD$ is real-analytic and compact there are at most 
finitely many points 
$\z\in bD$ such that the tangent line to $bD$ at $\z $ is parallel to 
the imaginary axis. 

For $z\in S$ write 
$$
\tilde E (z)=N\cap L(z),\ \ \tilde E_j (z)=N_j\cap L(z),\ \ j=1,2
$$
and 
$$
\Lambda (z) = \{ w\in\C\colon\ (z,w)\in N\},\ \ \Lambda _j (z)=
\{ w\in\C\colon\ (z,w)\in N_j\},\ j=1,2,  
$$
so that 
$$
\tilde E(z)=\{ z\} \times\Lambda (z),\ \ \tilde E_j (z) =
\{ z\}\times \Lambda _j (z),\ j=1,2.
$$ 
For each $t\in \R$ we have $N+t(i,0)=N,\ N_j + t(i,0) = N_j,\ j=1,2$,  so it follows that 
$$
\Lambda (z) =\Lambda (\Re z),\ \ \Lambda _j(z) =\Lambda _j (\Re z)\ \  (j=1,2,\ z\in S).
$$
Thus, it is enough to study $\Lambda (\tau ),\ \Lambda _j(\tau ),\ j=1,2$, where 
$\alpha\leq\tau\leq\beta$.

The set $\tilde E(\tau ) = \{\tau\}\times\Lambda (\tau )$, contained in 
$\R^3=\R\times\{ 0\}\times\C$, is the intersection of 
$\{ (\theta (w),w)\colon\ w\in\C\}$, the graph of 
$\theta$, with the two-plane (in fact, the complex line),
$\{\tau\}\times\C$. Since $\theta $ is
invariant with respect to the reflection across $b\D $ it follows that we get
$\Lambda _2(\tau)$ by reflecting 
$\Lambda _1 (\tau )\subset\overline\D$ across $b\D$. So it is enough 
to study $\Lambda _1 (\tau )$. 
Clearly 
$$\Lambda _1(\tau )=\Phi (K(\tau))\hbox{\ \ where\ \ }
K(\tau ) = (\tau +i\R )\cap\overline D.
$$
If $\tau \in [\alpha,\beta]$ is such that $\tau +i\R$ meets $bD$ 
transversely, as happens for all but finitely many $\tau$,
then $\Lambda _1(\tau )$ consists of finitely many pairwise disjoint 
closed arcs with endpoints on $b\D$ but otherwise contained in $\D $ which 
meet $b\D $ transversely. By transversality and by the fact that $\Phi $ 
extends across $bD$ as a conformal map, 
these arcs change smoothly with $\tau $ as long as $\tau +i\R $ meets $bD$ 
transversely. If
$\tau\in(\alpha,\beta )$ is such that $\tau +i\R $ does not meet  $bD$ 
transversely then $\Lambda (\tau )$ consists of a finite number of arcs with 
endpoints on $b\D $ and pairwise disjoint interiors plus a possible finite 
set on $b\D $. There may be points
on $b\D $ which are common endpoints of two (but not more than two ) of 
these arcs. Since we get $\Lambda _2(\tau )$ by
reflecting $\Lambda _1 (\tau )$ across $b\D $ it follows that if 
$\tau\not =0$ and if $\tau +i\R$ is transverse to $bD$ then $\Lambda (\tau )$ 
consists of finitely many pairwise disjoint simple closed curves, symmetric
with respect to 
$b\D $. If $\tau\not= 0$ and $\tau +i\R$ does not meet $bD$ transversely then 
$\Lambda (\tau )$ consists of finitely many pairwise disjoint simple closed
curves, symmetric 
with respect to $b\D$ plus a possible finite subset of $b\D $. There may 
be points on $b\D $
that are common points of two, 
but not more than two of these curves. Except 
for these points, the curves are pairwise disjoint. Clearly $\Lambda (\alpha )$ 
and $\Lambda (\beta )$ are finite sets. 

Since $i\R $ meets $bD$ transversely $\Lambda_1(0)$ consists of 
finitely many pairwise disjoint closed arcs with endpoints on $b\D $ but
otherwise contained in $\D $ which meet $b\D$ 
transversely. One of these arcs passes through the origin so its image
under the reflection across $b\D $ passes through infinity. Thus,
$\Lambda (0)\cup\{ \infty \}$ consists of finitely many pairwise 
disjoint simple closed curves on the 
Riemann sphere one of which contains infinity. 

For each $z,\ \alpha < \Re z<0$ the set $Y(z) = P_1\cap L(z)$ is a
bounded open
subset of $L(z)$ whose boundary $E(z) = bY(z)$ is the part of 
$\tilde E(z) = \{ z\}\times\Lambda (z)$ consisting of curves (recall 
that in addition 
to these curves, $\tilde E(z)$ may contain an additional finite set 
contained in $\{ z\} \times b\D )$. The complex line $L(z)$ has a natural 
orientation. 
We orient $E(z)$ as the boundary of $Y(z)$ in $L(z)$. Similarly, 
for $0<\Re z<\beta$ we orient $E(z)$ as the boundary of
$Y(z)=P_2\cap L(z)$ in $L(z)$. This determines the 
orientation of $\Lambda (z),\ 
\alpha < \Re z<\beta ,\ \Re z \not= 0$, or more precisely, the part
of $\Lambda (z)$ consisting of curves, and  
the orientation of $K(\tau ),\ \alpha <\tau <\beta , \tau\not= 0$, 
upwards if $\alpha<\tau<0$ and downwards if $0<\tau<\beta $.

If $0<\tau<\beta$ and a point $(\tau, w)$ is above the graph of $\theta $, 
that is, contained in $P_1$, then 
$[\tau,\beta]\times\{ w\}$ is contained in $P_1$. A consequence of this is
\vskip 2mm
\noindent\bf Proposition 5.1\ \it Suppose that $0<\tau<\beta $ and 
let $q_0\in L(\tau +
i\lambda (\tau))\cap P_1$. Then there is a path $t\mapsto q(t)\ (\tau\leq t\leq \beta)$
such that
$q(\tau )=q_0$ and $q(t)\in L(t+i\lambda (t))\cap P_1\ (\tau\leq t\leq \beta )$.\rm
\vskip 2mm
\noindent\bf Proof.\ \rm We have $q_0 = (\tau+i\lambda (\tau), w)$ 
where $\tau > \theta (w)$. Define $q(t)=(t+i\lambda (t),w)\ \ (\tau\leq t\leq \beta )$.
It is easy to see that 
$q$ has all the required properties.
\vskip 2mm
\noindent A similar proposition holds for $\alpha<\tau<0$ and for $P_2$ in the place of $P_1$. 
\vskip 4mm
6.\ \bf Continuity of an integral\rm
\vskip 2mm
We shall need
\vskip 2mm
\noindent\bf Proposition  6.1\it\ Let $z_0\in S$ and suppose that $G$ is a continuous
function on a neighbourhood of $\tilde E(z_0)$ in $N$. Then the function 
$$
\Theta (z)=\int _{\Lambda (z)} G(z,w) dw
$$
defined in a neighbourhood of $z_0$ in $S$, is continuous at $z_0$.
\vskip 2mm
\noindent\bf Proof.\ \rm We prove the continuity of 
$$
\Theta _1(z)=\int _{\Lambda _1(z)} G(z,w) dw =\int _{K(\Re z)}G(z,\Phi (\z))
\Phi^\prime (\z ) d\z .
$$
The proof for $\Theta _2(z)=\int _{\Lambda _2(z)} G(z,w) dw$ will be analogous; note that  
 $\Theta =\Theta _1+\Theta_2 $
since $\Lambda _1(z)$ and $\Lambda _2(z)$ meet in a finite set.
Recall that $\Phi $ and $\Phi^\prime $ extend holomorphically 
into a neighbourhood of $\overline D$, so the continuity of $\Theta _1$ depends on how
$K(t)=\overline D\cap (t+i\R )$  changes with $t$ near $t_0 =\Re z_0$. 
Assume for a moment that $\alpha < t_0<\beta $. There are at most finitely 
many $t,\ \alpha\leq t\leq \beta $ such that $t+i\R $ does not intersect $bD$ 
transversely. Thus there 
is an $\eta >0$ such that $t+i\R $ intersects $bD$ transversely 
for every $t,\ 0<|t-t_0|<\eta $. In particular, for each $t,\ t_0-\eta < t <\eta $, \ $K(t)$ 
is a finite collection of pairwise disjoint closed segments with 
endpoints varying continously with $t$. Since $bD\cap (t_0+i\R)$ is a finite 
set and since $bD$ is compact it  follows that each of these 
endpoints has a limit as $t\nearrow t_0$. As $t\nearrow t_0$ some segments
may degenerate into points in the limit $K(t_0-0)$ and some pairs of segments 
may get a common 
endpoint. Clearly $K(t_0-0)\subset K(t_0)$ and $K(t_0)\setminus 
K(t_0-0)\subset bD\cap (t_0+i\R)$. Since  the set $bD\cap (t_0+i\R)$ 
is finite it follows that $K(t_0-0)\subset K(t_0)$ is a finite set. Thus, 
$\lim_{z\rightarrow z_0, \Re z\leq\Re z_0} \Theta (z) = \Theta (z_0)$. 
Similarly we show that $\lim_{z\rightarrow z_0, \Re z\geq\Re z_0} 
\Theta (z) = \Theta (z_0)$ which proves that $\Theta $ is continuous at $z_0$. 
The same (one sided) reasoning applies if $z_0\in bS$. The proof is complete. 
 \vskip 4mm
 \bf 7.\ The manifold $M$ and the function $F$ \rm
 \vskip 2mm
 We now define a submanifold of $N$ that is more closely related to our problem. Write
 $$
 \eqalign{
 &M_1 = \{ (\xi +it, \Phi (\xi ))\colon\ \xi \in\overline D,\ a\leq t\leq b\} \cr
 &M_2 = \{ (\z +is,\ 1/\Psi (\z ))\colon\ \z\in \overline{D^\ast}\setminus \{ 0\},
 \ c\leq s\leq d\}\cr}
 $$
 and let $M=M_1\cup M_2$. Note that $M_1$ is a smooth manifold with boundary $bM_1=
 \{ (\xi +it,\Phi (\xi ))\colon\ \xi\in bD,\ a\leq t\leq b\}\cup 
 \{ (\xi +ia,\Phi (\xi ))\colon\ \xi\in\overline D\}\cup
\{ (\xi +ib,\Phi (\xi ))\colon\ \xi\in\overline D\}$. It is a submanifold of $N_1$.
Similarly, 
$M_2$ is a smooth manifold with boundary 
 $bM_2=
 \{ (\z +is, 1/\Psi (\z ))\colon\ \z\in bD^\ast ,\ c\leq s\leq d\}\cup 
 \{ (\z +ic, 1/\Psi (\z ))\colon\ \z\in\overline {D^\ast} \setminus\{ 0\} \}\cup
\{ (\z +id, 1/\Psi (\z ))\colon\ \z\in\overline {D^\ast} \setminus\{ 0\} \}.$

Suppose that $f$ is a continuous function 
on $Q_1\cup Q_2$ which satisfies (2.1) and (2.2). For each $t,\ a\leq t\leq b$, 
let $g_t$ be a continuous extension of 
$\xi\mapsto f(\xi +it)\ (\xi \in bD)$\ to $\overline D$ which is 
holomorphic on $D$ and for each $s,\ c\leq s\leq d,$ \ let $h_s$ be the continuous 
extension of $\z\mapsto f(\z + is)\ (\z\in bD^\ast )$ to $\overline{D^\ast}$ which 
is holomorphic on $D^\ast $. Define the function $G$ on $M_1$ by
$$
G(\xi +it,\Phi (\xi ))= g_t(\xi )\ \ (\xi\in\overline D,\ a\leq t\leq b)
$$
and the function $H$ on $M_2$ by
$$
H(\z+is,1/\Psi (\z ))= h_s (\z )\ \ (\z\in\overline{D^\ast}\setminus\{ 0\},\ c\leq s\leq d).
$$
In particular, on the part of $bM_1$ contained in $bN_1$ we have
$$
G(\xi +it, \Phi (\xi ))=f(\xi +it)\ \ (\xi\in bD,\ a\leq t\leq b),
$$
and on the part of $bM_2$ contained in $bN_2$ we have 
$$
H(\z +is,1/\Psi (\z ))= f(\z+is)\ \ (\z\in bD^\ast,\ c\leq s\leq d).
$$
Suppose that $(z,w)\in M_1\cap M_2$. Then there are $\xi\in bD,\ \z\in bD^\ast $ and $t, s,\ 
a\leq t\leq b,\ c\leq s\leq d$, such that
$(\xi +it,\Phi (\xi ))= (z,w)=(\z +is,1/\Psi (\z ))=
(\z+is, \overline{\Psi (\z)})=(\z +is,\Phi (\overline\z))$ which implies that
$\z=\overline\xi$ and
$\xi+it=\overline\xi +is$. Thus, $G(z,w)= G(\xi+it,\Phi (\xi))= f(\xi+it)=f(\overline\xi+is)=
H(\overline\xi+is, 1/\Psi (\overline \xi))= H(\z+is,1/\Psi (\z ))= H(z,w)$. It follows that
$$
F(z,w) = \left\{\eqalign{
&G(z,w)\ \ ((z,w)\in M_1) \cr
&H(z,w)\ \ ((z,w)\in M_2) \cr}\right.
$$
is a well defined continuous function on $M_1\cup M_2$ which
is holomorphic on each holomorphic leaf of $M_1$ and on each holomorphic leaf of $M_2$. 

Our aim is to show that $F$ depends only on $z$ which will imply that 
$f$ is holomorphic on $\hbox{Int}Q_1\cup \hbox{Int}Q_2$.
\vskip4mm
\noindent\bf 8.\ Integrals of CR functions on $M$\rm
\vskip 2mm
Denote by $\pi _1$ the projection $\pi_1(z,w)=z$. 
With no loss of generality assume that $\lambda (0)=0$. Our assumptions imply that 
there is an $\eta >0$ such that if
$$
\Sigma = \{ t+is\colon\ \lambda (t)-\eta<s<\lambda (t)+\eta,\ \alpha<t<\beta\} 
$$
then $L(z)\cap M =L(z)\cap N$ for all $z\in\Sigma $, that is, $\pi_1^{-1}(\Sigma)\cap M = 
\pi_1^{-1}(\Sigma)\cap N$.
Put $\Sigma _1 = \{ z\in\Sigma,\ \Re z<0\}$,\ $\Sigma _2 = \{ z\in\Sigma,\ \Re z>0\}$. 
Recall that for $z\in\Sigma _1$ the set $E(z)$ is the boundary of $Y(z)=P_1\cap L(z)$ in 
$L(z)$  and for $z\in\Sigma _2$ the set $E(z)$ is the boundary of $Y(z)=P_2\cap L(z)$ in 
$L(z)$. Let
$$
A_j=\cup\{ Y(z)\colon\ z\in\Sigma _j\} = \pi_1(\Sigma _j)\cap P_j\ \ (j=1, 2).
$$
For each $j=1, 2, \ A_j $ is an open 
subset of $\Sigma _j\times \C$ whose 
relative boundary is $N\cap (\Sigma_j\times\C)$. Using Proposition 5.1 
we see that the complement of 
$\overline A_j$ in $\Sigma _j\times \C$ is connected, $j=1, 2$. 

We shall prove that the function $F$ extends holomorphically into $A_1$ and into $A_2$.  
We begin to follow the reasoning of H.\ Lewy. In [L] this was done for smooth 
functions on smooth manifolds and for more general, including continuous, 
functions on smooth manifolds this was done in [R]. We cannot refer to these results directly 
 since in 
our case the manifold is not smooth but consists of two smooth pieces. However, these
two pieces are both foliated by analytic discs which simplifies the situation. 
We provide the details to make the proof self contained. 

Let $\D (u,r)$ be contained in either $\Sigma_1$ or $\Sigma _2$ and assume that 
$\Theta $ is a continuous function on $\pi _1^{-1}(\D (u,r))\cap N$
which is holomorphic on each holomorphic leaf. The function
 $$
 z\mapsto Q(z)=\int _{\Lambda (z)} \Theta (z,w)dw
 $$ 
 is, as we know, well defined and by Proposition 6.1 it is 
 continuous on $\D (u,r)$. Recall that there are 
 at most finitely many real
 values $\tau$  such that $\tau +i\R $ is not transversal to $bD$. So, 
 if we want to prove that $Q$ is holomorphic on $\D (u,r)$ it is enough to prove 
 that $Q$ is holomorphic in a neighbourhood of each $z_0\in \D (u,r)$ such that
 $z_0+i\R$ intersects  $bD$ transversely. Let $z_0$ be such a point. Let 
 $\rho>0$ be such that 
 for each $z\in \overline{\D (z_0,\rho )} $, $z+i\R$ meets $bD$ transversely. 
 Passing to a smaller $\rho $ 
 if necessary we may assume that there is a $\gamma >0$ such that whenever $U$ 
 is a closed disc contained in 
 $\D(z_0,\rho )$ of radius not exceeding $\gamma $, the circle $bU+it,\ t\in\R$, 
 either misses $bD$ 
 or meets $bD$ in one point or in two points. Let $U$ be such a disc. By 
 transversality, there is a positive 
 integer $\nu$ such that for each $z\in U$ the set $L(z)\cap N=L(z)\cap M$
 consists of $\nu $ pairwise disjoint 
 simple closed curves, each being the union of two smooth arcs with endpoints 
 on $\{ z\}\times b\D $, one having its interior contained in $\{ z\}\times \D$,
 and the other having its 
 interior contained in $\{ z\} \times (\C\setminus\overline\D )$, which change 
 smoothly with $z$. 
 So $\pi _1^{-1}(U)\cap N = \pi _1^{-1}(U)\cap M = \cup \{ L(z)\cap N\colon\ z\in U\} $ is 
 an open subset of $N$ which has $\nu $ components whose  closures are 
 pairwise disjoint;\ the boundary 
 of this set is 
 $\pi _1^{-1} (bU)\cap N = \cup \{ L(z)\cap N\colon z\in bU\}$. Let
 $\Omega$ be one of these components. Write $T = bN_1 = bN_2$. The set $\Omega 
 $ consists of three pairwise disjoint parts:\ the domains $\Omega_1=\Omega\cap\hbox{Int}N_1$, 
 $\Omega_2=\Omega\cap\hbox{Int}N_2$ and the two dimensional surface $\Omega\cap T$. For each 
 $z\in bU$, \ $L(z)\cap b\Omega $ is a simple closed curve so $b\Omega $ is a torus and $\Omega 
 $ is a solid torus in $N$. 
 
 \noindent We now want to show that
 $$
 \int_{b\Omega}\Theta (z,w) dz\wedge dw = 0.
 \eqno (8.1)
 $$
 Note first that $\Omega\cap T$ is the common part 
 of $b\Omega _1$ and $b\Omega _2$ so to prove (8.1) it is enough to 
 prove that 
 $$
  \int_{b\Omega_1}\Theta (z,w) dz\wedge dw = 0.
  \eqno (8.2)
 $$
and
$$
 \int_{b\Omega_2}\Theta (z,w) dz\wedge dw = 0.
 \eqno (8.3)
 $$
Consider (8.2). The properties of $U$ imply 
that $\Omega _1$ can be written as the union of a continous family of 
pairwise disjoint analytic discs
$$
A_t = \{ (\z+it,\Phi (\z ))\colon\ \z\in D, \z+it\in U\} =
\{((\z ,\Phi (\z ))+it\colon\ \z\in D\cap (-it+U)\}
$$
and $b\Omega_1$ is the union of their pairwise disjoint boundaries
$$
bA_t =\{ (\z, \Phi (\z ))+it\colon\z\in b(D\cap (-it+U))\}
$$
These analytic discs, if nonempty, are of two sorts: 
either their boundaries are smooth simple closed curves
which meet $\Omega\cap T$ in at most one point (which happens if $U\subset D+it$), 
or their boundaries are simple closed curves consisting of two arcs, one contained in 
$\hbox{Int} N_1$ and the other contained in $T$ (which happens if $U$ meets $D+it$ but is 
not contained in $D+it$). Recall that $N_1 = \{ (\Upsilon (w)+it, w)
\colon\ w\in\overline\D ,\ t\in\R\}$
where the conformal map $\Upsilon = \Phi^{-1} $ extends to a biholomorphic
map from $R\D $ for some $R>1$ 
to a neighbourhood of $\overline D$. Define the function $\rho $ by
$\rho (z,w)= (1/i)(z-\Upsilon (w))$ and notice that $\rho $ is real on $N_1$. 
Then $\Theta (z,w)dz\wedge dw= d\rho\wedge \mu$ 
where 
$\mu = -(1/i)(\Theta (z,w)/ \Upsilon^\prime (w)) dz$ which, by using the Fubini theorem on 
each of the three smooth pieces of $b\Omega _1$ and adding the results, implies that  
$$
\int_{b\Omega _1}\Theta (z,w)dz\wedge dw = \int_I\bigl[\int_{bA_t}\mu \bigr] dt  
\eqno (8.4)
$$
where $I\subset \R$ is the segment of all $t$ such that $\Omega_1\cap \{(\z+it,\Phi (\z ))\colon 
\ \z\in D\} 
$
is not empty. 
For each $t\in I$ we have 
$$
\int_{bA_t}\mu = -{1\over i}\int_{b A_t}\Theta (z,w){1\over{\Upsilon^\prime (w)}} dz
$$
where the integral on the right vanishes by the by the Cauchy theorem since the function
$(z,w)\mapsto \Theta (z,w)/\Upsilon^\prime (w)$ is 
continuous on $\overline A_t$ and holomorphic on 
$A_t$.  
 This proves that the integral on the left in (8.4) vanishes. 
We repeat the reasoning for $N_2$ to get (8.3). This proves (8.1). Thus, 
$$
\int_{bU}\Bigl[\int_{b\Lambda (z)} \Theta (z,w) dw\Bigr] dz = 0
$$
for every disc $U\subset \D (z_0,\rho )$, which, by the Morera theorem, implies that the 
function $Q$ is holomorphic on $D(z_0,\rho )$. This proves that $Q$ is 
holomorphic on $\Sigma_1\cup\Sigma _2$. 
\vskip 4mm
\bf 9.\ Holomorphic extensions of $F$ and the completion of the proof\rm 
\vskip 2mm
We continue to follow the reasoning of H.\ Lewy. For each $z\in\Sigma_1\cup\Sigma_2$ 
and each $w\in\C\setminus 
\Lambda (z)$ (that is, for each $w$ such that $(z,w)\not\in N)$ define 
$$
\Xi (z,w) = {1\over{2\pi i}}\int _{\Lambda (z)}{{F(z,\z )}\over{\z - w}} d\z .
$$
For a fixed $z$ the function $w\mapsto \Xi (z,w)$ is holomorphic on $\C\setminus \Lambda (z)$. 
Now, fix $z_0\in\Sigma_1$
and $w\in\C\setminus \Lambda (z_0)$. We show that $z\mapsto \Xi (z,w)$ is 
holomorphic in a neighbourhood of
$z_0$. There is a $\rho >0$ such that $\D(z_0,\rho ) \subset\Sigma_1$ and $(z,w)\not\in N\ (z\in 
\D (z_0,\rho ))$ so the function $(z,\z )\mapsto F(z,\z )/(\z -w)$ is continuous on 
$\cup \{L(z)\cap N\colon\ z\in \D (z_0,\rho )\}$ and holomorphic on each holomorphic leaf. By 
the 
reasoning in Section 8 it follows that $z\mapsto \Xi (z, w )$ is holomorphic on $\D (z_0,\rho)$. 
This shows that $\Xi $ is holomorphic on $\pi_1^{-1}(\Sigma_1)\setminus N$. Fix a large $w$.
We know that 
 $z\mapsto \Xi (z,w)$ is continuous on $\Sigma _1\cup [\alpha + i(-\delta ,\delta)]$ 
and holomorphic on $\Sigma $. Since $\Lambda (\alpha )$ is a finite set it follows that $
\Xi (z,w)$ approaches $0$ as $z$ approaches a point of $\alpha+i(-\delta,\delta)$. It follows 
that
$\Xi (z,w)\equiv 0\ (z \in \Sigma _1)$.  Since this holds for all sufficiently large $w$,
the connectedness of
$(\Sigma_1\times \C)\setminus\overline{A_1}$ implies that 
$\Xi \equiv 0$ on $(\Sigma_1\times \C)\setminus\overline{A_1}$ so
$$
{1\over{2\pi i}}\int_{\Lambda (z)} {F(z,\z )\over {\z - w}} d\z \equiv 0\ \ 
(z\in\Sigma_1,\ w\in \C\setminus \overline{Y(z)} ).
$$
It follows, in particular, that for all $z\in\Sigma_1$ the function $w\mapsto F(z,w)$, 
defined on $E(z)=bY(z)$, has a continuous extension into $Y(z)\cup bY(z)$ which is 
holomorphic on $Y(z)$. The same reasoning shows this for all $z\in \Sigma _2$. 

Recall that $\Lambda (0)\cup\{ \infty \}$ consists of finitely many pairwise 
disjoint simple closed curves on the 
Riemann sphere one of which passes through infinity. 
By transversality, $\Lambda (t)$ changes continuously with $t$ near 
$0$ and contains $\infty $ only if $t=0$. As 
$t\nearrow 0$ the open sets $Y(t)$ and their oriented boundaries 
converge to a domain $Y^-$ and 
its oriented boundary $bY^-$ which, as a set, coincides with 
$\Lambda (0)$. As $t\searrow 0$ the 
open sets $Y(t)$ and their oriented boundaries $\Lambda (t)$ 
converge to a domain $Y^+$ and its 
oriented boundary $bY^+$ which, as a set, coincides with $\Lambda (0)$. 
We have $Y^-\cup\Lambda(0)\cup Y^+ =\C$. Since the 
function $F$ is bounded and continuous on $M$ it follows that as $t \nearrow 0$, 
the holomorphic extensions of
$F|M\cap L(t+i\lambda (t))$ converge to the holomorphic extension of $F|M\cap L(0)$, 
a bounded continuous function on 
$\{ 0\}\times\overline{Y^-} $, holomorphic on $\{ 0\} \times Y^-$. In particular, 
$F|M\cap L(0)$ extends to a bounded continuous function on $\{ 0\}\times\overline{Y^-} $ 
which is holomorphic on 
$\{ 0\}\times {Y^-} $. Similarly, $F|M\cap L(0)$ extends to a bounded continuous function 
on 
$\{ 0\}\times\overline{Y^+} $ which is holomorphic on $\{ 0\}\times\overline{Y^+} $. 
Consequently $F|M\cap L(0) $ extends to 
a bounded continuous function on $\{ 0\}\times \C$ 
which is holomorphic  on $\{ 0\} \times [Y^+\cup Y^-]$. This function 
extends holomorphically across $\{ 0\}\times \Lambda (0)$ to a bounded entire function
on $L(0)$, which, 
by the Liouville theorem, must be constant. This implies that $F|M\cap L(0)$ is a 
constant. In the same way we show that $F|M\cap L(is)$ is constant for each 
$s, \ -\eta<s<\eta$. It follows that there is an
$r>0$ such that on $r\D $ the
holomorphic extensions of all functions $f|(it+bD)$ to $it+D$, for all $t$ 
such that $r\D\subset it+D$, coincide, and that 
on $r\D$  the
holomorphic extensions of all functions $f|(it+bD^\ast )$ to $it+D^\ast$, for all $t$ 
such that $r\D\subset it+D^\ast$, coincide. This implies first that $f$ is
holomorphic on $\cup\{it+bD
\colon r\D\subset it+D,
\ a<t<b\}\cap \hbox{Int}S$, 
and then, by translating $bD$ further along $i\R$, that 
$f$ is holomorphic on $\cup\{ it+bD\colon \ a<t<b\} 
\cap \hbox{Int}S =\hbox{Int}Q_1$.
Similarly, we show that $f$ is holomorphic on 
$\hbox{Int}Q_2$.  The proof is complete. 
\vskip 4mm
\bf 10.\ Remarks \rm
\vskip 2mm
It remains an open problem to prove Theorem 2.1 without the assumption
that $\Omega $ is symmetric with respect to the real axis. 
It is known that one cannot drop the assumption 
that $(\overline\Omega -ir)\cap(\overline\Omega +ir)=
\emptyset $. In fact the function
$$
f(z)=\left\{\eqalign{&z^2/\overline z\ \ (z\not= 0)\cr
&0\ \ \ \ \ \ \ (z=0)\cr}\right.
$$ 
is continuous on $\C$ and extends holomorphically from 
each circle which either surrounds the origin or contains the origin, yet $f$ is not holomorphic.
\vskip 10mm
\bf \noindent Acknowledgements\ \ \rm The author is indebted to Mark Agranovsky for several 
very sti\-mu\-la\-ting discussions. After these discussions the author found a way to pass from 
circles to general curves. 

This work was supported 
in part by the Ministry of Higher Education, Science and Technology of Slovenia 
through the research program Analysis and Geometry, Contract No.\ P1-0291
\vfill
\eject
\centerline{\bf References}
\vskip 2mm
\noindent [AG]\ M.\ L.\ Agranovsky and J.\ Globevnik: Analyticity 
on circles for rational and real-analytic functions of two real variables.

\noindent J.\ d'Analyse Math.\ 91 (2003) 31-65
\vskip 2mm
\noindent [E]\ L.\ Ehrenpreis: Three problems at Mount Holyoke. 

\noindent Contemp.\ Math.\ 278 (2001) 123-130. 
\vskip 2mm
\noindent [G1] J.\ Globevnik:\ Analyticity on rotation invariant families of curves.

\noindent Trans.\ Amer.\ Math.\ Soc.\ 280 (1983) 247-257
\vskip 2mm
\noindent [G2] J.\ Globevnik: Analyticity on familes of curves.

\noindent Talk at Bar Ilan University, November 1987
\vskip 2mm
\noindent [G3] \ J.\ Globevnik: Holomorphic extensions from open families of circles.

\noindent Trans.\ Amer.\ Math.\ Soc.\ 355 (2003) 1921-1931
\vskip 2mm
\noindent [L]\ H.\ Lewy: On the local character of the 
solutions of an atypical linear differential equation in three 
variables and a related theorem for regular functions of two complex variables.

\noindent Ann. of Math.\ 64 (1956) 514-522
\vskip 2mm
\noindent [R]\ H.\ Rossi: A generalization of a theorem of Hans Lewy.

\noindent Proc.\ Amer.\ Math.\ Soc.\ 19 (1968) 436-440
\vskip 2mm
\noindent [T1]\ A.\ Tumanov: A Morera type theorem in the strip.

\noindent Math.\ Res.\ Lett.\ 11 (2004) 23-29
\vskip 2mm
\noindent [T2]\ A.\ Tumanov: Testing analyticity on circles.

\noindent \ Preprint [ArXiv:math.CV/0502139]
\vskip 10mm
\noindent Institute of Mathematics, Physics and Mechanics

\noindent University of Ljubljana, Ljubljana, Slovenia

\noindent e-mail: josip.globevnik@fmf.uni-lj.si

\bye